\documentclass[a4paper,leqno,11pt,final]{amsproc}
\usepackage{amssymb,latexsym,amsxtra,amscd,amsfonts,amsthm,amsmath,verbatim}

\newcommand{\lm}{{\rm Lm}}

\numberwithin{equation}{section}

\theoremstyle{remark}

\newtheorem{theorem}[equation]{Theorem}
\newtheorem{proposition}[equation]{Proposition}
\newtheorem{lemma}[equation]{Lemma}
\newtheorem{corollary}[equation]{Corollary}

\newtheorem{remark}[equation]{Remark}

\begin{document}

\title{THE REES ALGEBRA FOR CERTAIN MONOMIAL CURVES}

\author{Debasish  Mukhopadhyay}
\address{Acharya Girish Chandra Bose College, 35, Scott Lane, Kolkata, WB 700009, INDIA.}
\email{mdebasish01@yahoo.co.in}
\author{Indranath Sengupta}
\thanks{Research supported by the DST Project No. SR/S4/MS: 614/09}
\address{School of Mathematical Sciences, Ramakrishna Mission Vivekananda University, Belur Math, Howrah, WB 711 202, INDIA}
\curraddr{DEPARTMENT OF MATHEMATICS, JADAVPUR UNIVERSITY, KOLKATA,
WB 700 032, INDIA.}
\email{sengupta.indranath@gmail.com}

\begin{abstract}
Let $K$ be a perfect field and let $m_{0} < m_{1} < m_{2} < m_{3}$ be a
sequence of coprime positive integers such that they form a minimal
arithmetic progression. Let $\wp$ denote the defining ideal of the
monomial curve $\mathcal{C}$ in $\mathbb{A}_{K}^{4}$, defined by the
parametrization $X_{0} = T^{m_{0}}, X_{1} = T^{m_{1}},X_{2} = T^{m_{2}},X_{3} =
T^{m_{3}}$. Let $R = K[X_{0}, X_{1}, X_{2}, X_{3}]$. In this article,
we find the equations defining the Rees
algebra $R[\wp t]$ explicitly and use them to prove that the blowup
scheme ${\rm Proj}\, R[\wp t]$ is not smooth. This proves Francia's
conjecture in affirmative, which says that a dimension one prime in
a regular local ring is a complete intersection if it has a smooth
blowup.

\begin{center}
Keywords~:~ Monomial Curves\,,\, Gr\"{o}bner Basis\,,\, Rees Algebra\,.\\

Mathematics Subject Classification 2000~:~ 13P10\,,\,13A30\,. \\
\end{center}

\end{abstract}

\maketitle

\section{Introduction}
\noindent Blowup Algebras, in particular the Rees algebra $\mathcal{R}(I) = R[It]$
($t$ a variable) and the associated graded ring $\mathcal{G}(I) =
\mathcal{R}(I)/I\mathcal{R}(I)$ of an ideal $I$ in a Noetherian ring
$R$ play a crucial role in the birational study of algebraic
varieties. The scheme ${\rm Proj}(\mathcal{R}(I))$ is the blowup of
${\rm Spec}(R)$ along $V(I)$, with ${\rm Proj}(\mathcal{G}(I))$ being
the exceptional fiber. Although blowing up is a
fundamental operation, an explicit understanding of this process
remains an open problem. For example, Francia's conjecture stated in
O'Carroll-Valla (1997) says: If $R$ is a regular local ring and $I$
is a dimension one prime ideal in $R$ then $I$ is a complete
intersection if ${\rm Proj}(\mathcal{R}(I))$ is a smooth projective
scheme. A negative answer to this conjecture was given by
Johnson-Morey (2001; 1.1.5), for $R = \mathbb{Q}[x,y,z]$. It is
still unknown whether the conjecture is true or not for polynomial
rings $R$ over an algebraically closed field. It is evident that a
good understanding of the defining equations of the Rees algebra is
necessary to answer such queries and an explicit computation of
these equations is often extremely difficult. In this context,
the Cohen-Macaulay and the normal properties of blowup algebras
have attracted the attention of several authors because they help
in describing these algebras qualitatively.
\medskip

Our aim in this article is to use the Elimination theorem to explicitly
compute the defining equations of the Rees algebra for certain one dimensional
prime ideals $\wp$, namely those which arise as the defining ideal of the affine
monomial curve given by the parametrization $X_{0} = T^{m_{0}}, X_{1} = T^{m_{1}},
X_{2} = T^{m_{2}}, X_{3} = T^{m_{3}}$ such that
$m_{0} < m_{1} < m_{2} < m_{3}$ is a sequence of positive
integers with gcd $1$, which form an arithmetic progression (see Section 4 for
complete technical details). The explicit form of these
equations will be used in Section 6, in conjuction with the Jacobian Criterion for smoothness
over a perfect field to prove that ${\rm Proj}R[\wp t]$ is not smooth. It
is known from the work of Maloo-Sengupta (2003), that $\wp$ is not a
complete intersection. Hence Francia's conjecture is true for $\wp$
over any perfect field $K$.

\section{Equations defining the Rees Algebra}

\noindent In order to compute the equations defining the Rees
algebra $\mathcal{R}(I)$, we view $\mathcal{R}(I)$ as quotients of
polynomial algebras. Thus for a Rees Algebra $\mathcal{R}(I)$, it amounts to
the study of the natural homomorphism associated to the
generators $(a_{1},\ldots, a_{m})$ of $I$
$$\widehat{R}=R[T_{1},\ldots,T_{m}]\stackrel{\varphi}\longrightarrow R[It], \quad \varphi(T_{i}) = a_{i}t\,;$$
and particularly of how to find $E = ker(\varphi)$, and analyze its
properties. $E$ will be referred to as the {\it equations of}
$\mathcal{R}(I)$ or {\it the defining ideal} of $\mathcal{R}(I)$.
One approach to get at these equations goes as follows. Let\\
$$ R^{r}\stackrel{\varphi}\longrightarrow
R^{m}\longrightarrow I \longrightarrow 0\,,$$ be a presentation of
the ideal $I$. $E_{1}$ is generated by the $1$-forms
$$[f_{1},\ldots,f_{r}] = [T_{1},\ldots,T_{m}] . \varphi = \mathbf{T}.\varphi.$$
The ring $R[T_{1},\ldots,T_{m}] / (E_{1})$ is the symmetric algebra
of the ideal $I$, and we write $\mathcal{A} = E / (E_{1})$ for the
kernel of the canonical surjection
$$0 \longrightarrow \mathcal{A} \longrightarrow S(I) \longrightarrow \mathcal{R}(I) \longrightarrow
0.$$ If $R$ is an integral domain, $\mathcal{A}$ is the $R$-torsion
submodule of $S(I)$. The ideal $I$ is said to be an {\it ideal of
linear type} if $\mathcal{A} = 0$, i.e., $E = E_{1}$, or
equivalently the symmetic algebra and the Rees algebra are
isomorphic. We will come accross a natural class of such ideals in
Section 6. An ideal reference is Vasconcelos (1994) for more on Blowup algebras and
related things.

\section{Computational Methods}

In this section, we assume the basic knowledge of Gr\"{o}bner bases and recall the
Elimination Theorem below, mostly from the book Cox-Little-O'Shea (1996).

\subsection{The Elimination Theorem}

Let $K[t_{1}, \ldots, t_{r}, Y_{1}, \ldots, Y_{s}]$ be a polynomial
ring over a field $K$. Let $\mathfrak{a}$ be an ideal in $K[t_{1},
\ldots, t_{r}, Y_{1}, \ldots, Y_{s}]$. The $r$-th elimination ideal
is $\mathfrak{a}(r) = \mathfrak{a} \cap K[Y_{1}, \ldots, Y_{s}]$.
We can actually compute a Gr\"{o}bner basis for $\mathfrak{a}(r)$,
if we know that of $\mathfrak{a}$ and if we choose a monomial order
suitably on $K[t_{1}, \ldots, t_{r}, Y_{1}, \ldots, Y_{s}]$.
Let $>_{\mathcal{E}}$ be a monomial order on $K[t_{1},
\ldots, t_{r}, Y_{1}, \ldots, Y_{s}]$, such that
$$t_{1} >_{\mathcal{E}} \cdots >_{\mathcal{E}} t_{r}
>_{\mathcal{E}} Y_{1} >_{\mathcal{E}}\cdots >_{\mathcal{E}} Y_{s}$$
and monomials involving at least one of the $t_{1}, \ldots, t_{r}$
are greater than all monomials involving only the remaining
variables $Y_{1}, \ldots, Y_{s}$. We then call $>_{\mathcal{E}}$ an
{\it elimination order} with respect to the variables $t_{1},
\ldots, t_{r}$.

\noindent One of the main tools for computing the equations of the
Rees algebras is the {\it Elimination Theorem}, which is the
following:

\begin{theorem}{\it
Let $G$ be a Gr\"{o}bner basis for the ideal $\mathfrak{a}$ in $K[t_{1}, \ldots, t_{r},
Y_{1}, \ldots, Y_{s}]$, where the order is an
elimination order $>_{\mathcal{E}}$ with respect to the variables $t_{1}, \ldots,
t_{r}$. Then $G_{r} = G \,\cap\, K[Y_{1}, \ldots, Y_{s}]$ is a
Gr\"{o}bner basis of the $r$-th elimination ideal
$\mathfrak{a}(r)$, with respect to $>_{\mathcal{E}}$. }
\end{theorem}

\proof See Cox-Little-O'Shea (1996; Chapter 3).\qed
\medskip

\noindent Let $I = (a_{1}, \ldots, a_{m})$ be an ideal in the
polynomial ring $R := K[Z_{1},\ldots, Z_{n}]$, over a field $K$. The
presentation of the Rees algebra $R[It]$ is obtained as:

\begin{proposition}{\it
In the ring $R[z_{1}, \ldots, z_{m}, t]$, consider the ideal
$\mathfrak{a}$ generated by the polynomials $z_{j} - ta_{j}$, $j =
1, \ldots, m$. Then $R[It] = R[z_{1}, \ldots, z_{m}]/E$, where $E =
\mathfrak{a} \cap R[z_{1}, \ldots, z_{m}]$. }
\end{proposition}

\proof It is clear that $E \supset \mathfrak{a} \cap R[z_{1},
\ldots, z_{m}]$. Conversely, if $f(z_{1}, \ldots, z_{m})$ is an
element of $E$, we write
$$f(z_{1}, \ldots, z_{m}) = f\left(ta_{1} + (z_{1} - ta_{1}), \ldots, ta_{m} + (z_{m} - ta_{m})\right)$$
and we can use Taylor expansion to show that $f\in
\mathfrak{a}$.\qed

\begin{proposition}{\it Let $R$, $\mathfrak{a}$ and $E$ be as defined in the Proposition 1.2.3.
Let $>_{\mathcal E}$ be an elimination order with respect to the
variable $t$ on $R[z_{1}, \ldots, z_{m}, t]$, with $t
>_{\mathcal{E}} \,Z_{i}, \,z_{j}$. If $\mathcal{G}$ is a Gr\"{o}bner
basis for $\mathfrak{a}$ with respect to $>_{\mathcal{E}}$, then
$\mathcal{G}\cap R[z_{1}, \ldots, z_{m}] $ is a Gr\"{o}bner basis
for $E$. }
\end{proposition}

\proof Follows from Theorem 3.1 and Proposition 3.2. \qed

\noindent We end this section with the statement of the Jacobian
Criterion for smoothness, which will be used for verifying
smoothness of the blowup; see Kunz (1985; page 171) for a proof.
\begin{theorem}{\it
Let $R=K[Z_{1},\ldots, Z_{n}]$ be a polynomial ring over a perfect field
$K$. Let $J=(f_{1},\ldots, f_{m})$ be an ideal in $R$ and set $S=R/J$. Let
$\mathfrak{p}$ be a prime ideal of $R$ containing $J$ and write
$\kappa(\mathfrak{p})=K(R/\mathfrak{p})$ for the residue field at
$\mathfrak{p}$. Let $c$ be the codimension of $J_{\mathfrak{p}}$ in
$R_{\mathfrak{p}}$.
\begin{enumerate}

\item The Jacobian matrix
$$\mathcal{J} := (\partial f_{i}/\partial Z_{j}),$$
taken modulo $\mathfrak{p}$ has rank atmost $c$.

\item $S_{\mathfrak{p}}$ is a regular local ring iff the matrix
$\mathcal{J}$, taken modulo $\mathfrak{p}$, has rank $c$.

\end{enumerate}
}
\end{theorem}

\section{Monomial Curves}

\noindent Let $\mathbb{N}$ \, and $\mathbb{Z}$ denote the set of
nonnegative integers and the set of integers respectively.
Assume that $0<m_{0} < m_{1} < \ldots < m_{p} $ form an
arithmetic sequence of integers, with $p\geq 2$ and $\gcd(m_{0},
\ldots ,m_{p}) = 1$. We further assume that $m_{i}=m_{0}+id$ where
$d$ is the common difference of the arithmetic sequence and $m_{0},
m_{1}, \ldots , m_{p}$ generate the numerical semigroup $\Gamma :=
\sum_{i=0}^{p}\mathbb{N}m_{i}$ minimally . Write $m_{0}=ap+b$, where
a and b are unique integers such that $a \geq 1$ (otherwise $m_{0},
m_{1}, \ldots , m_{p}$ can not generate the numerical semigroup $\Gamma$
minimally) and $1\leq b \leq p$. Let $\wp$ denote the kernel of the map $\eta : R:= K[X_{0}, X_{1}, \ldots,
X_{p}] \to K[T]$, given by $\eta(X_{i}) = T^{m_{i}} $. The
prime ideal $\wp$ is an one-dimensional perfect ideal and it is the
defining ideal of the affine monomial curve given by the
parametrization $X_{0} = T^{m_{0}}, \ldots, X_{p} = T^{m_{p}}$.
A minimal binomial generating set $\mathcal{G}$ for $\wp$
was constructed by Patil (1993). It was proved by Sengupta (2003) that
it is a Gr\"{o}bner basis with respect to the graded reverse lexocographic
monomial order. It was noted in Maloo-Sengupta (2003) that the set
$\mathcal{G}$ depends intrinsically on the integer $b$. We therefore
write $\mathcal{G}_{b}$ instead of $\mathcal{G}$, which is
$\mathcal{G}_{b}:=\{\phi(i,j) \mid i,j \in [1,p-1]\} \cup \{\,\psi(b, \,j) \mid j\in [0,p-b]\}$, such that\footnote{Our notations differ slightly from those introduced by Patil (1993) in the following manner: The embedding dimension in our case is $p+1$ and not $e$; the indeterminates $X_{0}, \ldots, X_{p}, Y$ have been replaced by $X_{0}, \ldots, X_{p}$; the binomials $\xi_{ij}$ occur in our list of binomials $\phi(i,j)$; the binomial $\theta$ is $\psi(b,p-b)$ in our list.}:

\begin{enumerate}
\item[(i)]   $\phi(i,j):=\begin{cases} X_{i}X_{j}-X_{\epsilon(i,j)}X_{i+j-\epsilon(i,j)}\,, & {\rm if} \quad i,j \in [1,p-1];\\
0\,, & {\rm otherwise}\,;
\end{cases}$
\medskip

\item[(ii)]  $\psi(b,j):=\begin{cases} X_{b+j}X_{p}^{a}-X_{j}X_{0}^{a+d}\,, & {\rm if} \quad j \in [0,p-b];\\
0\,, & {\rm otherwise}\,;
\end{cases}$
\end{enumerate}
\medskip

with
\medskip

\begin{enumerate}
\item[(iii)] $\epsilon(i\,,\,j) :=\quad
\begin{cases}
i+j & {\rm if} \quad i+j \,<\, p\\
p & {\rm if} \quad i+j \,\geq\, p\\
\end{cases}$;
\medskip

\item[(iv)] $[a\,,\,b]=\{i\in\mathbb{Z} \,\mid \,a\leq i\leq b\}$.
\end{enumerate}
\medskip

We now restrict our attention to $p=3$, since we will be
dealing only with monomial curves in affine
$4$-space, parametrized by four integers $m_{0}, \ldots , m_{3}$ in
arithmetic progression.
Let us write, $R_{b} = K[\mathbb{X}, \Psi_{b},
\Phi],$ such that $\Psi_{b}= \{\Psi(b,0),
\Psi(b,1),\ldots,\Psi(b,3-b)\}$, $\Phi=
\{\Phi(2,2),\Phi(1,2),\Phi(1,1)\}$ and
$\mathbb{X}=\{X_{1}, X_{2}, X_{3}, X_{0}\}$ are indeterminates. The indeterminate $X_{0}$ in the
set $\mathbb{X}$ has been listed at the end deliberately, keeping the monomial order in mind, to be
defined in the next section.
\medskip

Let \,$t$\, be an
indeterminate. We define the homomorphism $\varphi_{b}:
R_{b}\longrightarrow R[\wp t]$ as $\varphi_{b}(X_{i}) = X_{i}$, $\varphi_{b}(\Phi(i,j)) = \phi(i,j)t$,
$\varphi_{b}(\Psi(b,j))=\psi(b,j)t$. Let $E_{b}$ denote the kernel of
$\varphi_{b}$. Our aim is to construct a minimal Gr\"{o}bner basis
for the ideal $E_{b}$. Write $S = R_{b}[t]$ and define the ring
homomorphism \,$\overline{\varphi_{b}} : S \longrightarrow R[\wp t]$ as $\overline{\varphi_{b}}(t)=t \quad {\rm and} \quad
\overline{\varphi_{b}}=\varphi_{b} \quad {\rm on} \quad
R_{b}$. We follow the method of elimination described in
Propositions 3.2 and 3.3 and consider the ideal $\mathfrak{a}_{b}\subseteq S$ such that $\mathfrak{a}_{b} \cap R_{b} = E_{b}$.
We shall compute a Gr\"{o}bner basis $\widehat{\mathfrak{a}_{b}}$ for
$\mathfrak{a}_{b}$, with respect to an elimination order
$>_{\mathcal{E}}$ (with respect to $t$) on $S$. Then,
$\widehat{\mathfrak{a}_{b}}\cap R_{b}$ is a Gr\"{o}bner basis
for $E_{b}$, that is, those elements of $\widehat{\mathfrak{a}_{b}}$
that do not involve the variable $t$. These generators of $E_{b}$ will be used
to decide the non-smoothness of the blowup in section 6. We now define the desired
elimination order on $S$.

\section{Elimination order on $S = R_{b}[t]$}
\noindent A monomial in $\displaystyle{S =
R_{b}[t]=R[\Psi(b,0)
,\ldots,\Psi(b,3-b),\Phi(2,2),\Phi(1,2),\Phi(1,1),t]}$ is given by
$$\displaystyle{t^{d}\mathbb{X}^{\alpha}\Psi^{\beta}\Phi^{\gamma}
=
t^{d}\left(X_{1}^{\alpha_{1}}X_{2}^{\alpha_{2}}X_{3}^{\alpha_{2}}X_{0}^{\alpha_{0}}\right)
\left(\prod_{i=0}^{3-b}\Psi(b,i)^{\beta_{i}}\right)
\left(\Phi(2,2)^{\gamma_{1}}\Phi(1,2)^{\gamma_{2}}\Phi(1,1)^{\gamma_{3}}\right)
}$$ which is being identified with the ordered tuple
$\displaystyle{(d , \alpha, \beta,  \gamma)\in \mathbb{N}^{12 -b}}$,
such that
$$\displaystyle{\alpha := (\,\alpha_{1}, \alpha_{2},
\alpha_{3}, \alpha_{0}\,)}, \displaystyle{\quad \beta := (\,\beta_{0},\ldots, \beta_{3-b}\,)},
\displaystyle{\quad \gamma
:=(\,\gamma_{1},\gamma_{2},\gamma_{3}\,)}.$$.

\noindent Let us define the weight function
$\displaystyle{\widehat{\omega}}$ on the non-zero monomials of
$\displaystyle{S}$ to be the function with the property $\displaystyle{\widehat{\omega}(fg)=\widehat{\omega}(f)+\widehat{\omega}(g)},$
for any two non-zero monomials $f$ and $g$ in $S$, and that
$$\displaystyle{\widehat{\omega}(t)\,=\,1}, \quad \displaystyle{\widehat{\omega}(X_{i})\,=\,m_{i}}, \quad \displaystyle{\widehat{\omega}(\Phi(i,j))\,=\,\widehat{\omega}(X_{i}X_{j})}, \quad \displaystyle{\widehat{\omega}(\Psi(b,j))\,=\,\widehat{\omega}(X_{3}^{a}X_{b+j})}.$$

\noindent We say that $\quad\displaystyle{t^{d}\mathbb{X}^{\alpha}\Psi^{\beta}\Phi^{\gamma} \,
>_{\mathcal{E}} \,
t^{d'}\mathbb{X}^{\alpha'}\Psi^{\beta'}\Phi^{\gamma'} }, \quad$ if one of the following holds:

\begin{enumerate}
\item[(i)] $\displaystyle{d>d'}$;
\medskip

\item[(ii)]  $\displaystyle{d=d'}$ \,and\, $\displaystyle{\widehat{\omega}(\mathbb{X}^{\alpha}\Psi^{\beta}\Phi^{\gamma})
\,>\,\widehat{\omega}(\mathbb{X}^{\alpha}\Psi^{\beta}\Phi^{\gamma})}$;
\medskip

\item[(iii)]  $d=d'$, \,$\displaystyle{\widehat{\omega}(\mathbb{X}^{\alpha}\Psi^{\beta}\Phi^{\gamma})
\,=\,\widehat{\omega}(\mathbb{X}^{\alpha}\Psi^{\beta}\Phi^{\gamma})}$ \,and\, $\displaystyle{\sum\beta_{i}\,>\,\sum\beta'_{i}}$;
\medskip

\item[(iv)]  $\displaystyle{d=d'}$, \,$\displaystyle{\widehat{\omega}(\mathbb{X}^{\alpha}\Psi^{\beta}\Phi^{\gamma})
\,=\,\widehat{\omega}(\mathbb{X}^{\alpha}\Psi^{\beta}\Phi^{\gamma})}$, \,$\displaystyle{\sum\beta_{i}\,=\,\sum\beta'_{i}}$
\,and in the difference \,$\displaystyle{( \beta-\beta')}$, the rightmost non-zero entry is negative;
\medskip

\item[(v)]  $d=d'$, \,$\displaystyle{\widehat{\omega}(\mathbb{X}^{\alpha}\Psi^{\beta}\Phi^{\gamma})
\,=\,\widehat{\omega}(\mathbb{X}^{\alpha}\Psi^{\beta}\Phi^{\gamma})}$, \,$\displaystyle{ \beta=\beta'}$ \,and\, $\displaystyle{\sum\gamma_{i}\,>\,\sum\gamma'_{i}}$;
\medskip

\item[(vi)]  $d=d'$, \,$\displaystyle{\widehat{\omega}(\mathbb{X}^{\alpha}\Psi^{\beta}\Phi^{\gamma})
\,=\,\widehat{\omega}(\mathbb{X}^{\alpha}\Psi^{\beta}\Phi^{\gamma})}$, \,$\displaystyle{ \beta=\beta'}$, \,$\displaystyle{\sum\gamma_{i}\,=\,\sum\gamma'_{i}}$ \,and in the difference\, $\displaystyle{( \gamma-\gamma')}$, the rightmost
non-zero entry is negative;
\medskip

\item[(vii)]  $d=d'$, \,$\displaystyle{\widehat{\omega}(\mathbb{X}^{\alpha}\Psi^{\beta}\Phi^{\gamma})
\,=\,\widehat{\omega}(\mathbb{X}^{\alpha}\Psi^{\beta}\Phi^{\gamma})}$, $\displaystyle{\beta=\beta'}$,
\,$\displaystyle{\gamma=\gamma'}$ \, and in the difference \,$\displaystyle{(\alpha - \alpha')}$,
\,the rightmost non-zero entry is negative.
\end{enumerate}

\noindent  Then $>_{\mathcal{E}}$ is the desired elimination order on $S$, with respect to the variable $t$.
\medskip

\begin{theorem}
{\it Given $b\in \{1, 2,3\}$, let $\mathfrak{a}_{b}$ be the ideal in
$S$, generated by

\begin{itemize}
\item $P(i,j)=\begin{cases}\underline{tX_{i}X_{j}}
-tX_{\epsilon(i,j)}X_{i+j-\epsilon(i,j)}
-\Phi(i,j)\quad , \quad i,j \,\in \,[1,2]\,,\\[2mm]
0\hspace*{2.50in},\hspace*{0.20in} {\rm otherwise}\,;
\end{cases}$\\[3mm]

\medskip

\item $P(\Psi(b,l))=\begin{cases}\underline{tX_{b+l}X_{3}^{a}}-tX_{l}X_{0}^{a+d}-\Psi(b,l)
\quad \quad , \quad l\in
[0,3-b]\,,\\[2mm]
0\hspace*{2.30in},\hspace*{0.20in} {\rm otherwise}\,;
\end{cases}$
\end{itemize}

\noindent A Gr\"{o}bner Basis for the ideal $\mathfrak{a}_{b}$ is the set
$$\widehat{\mathfrak{a}_{b}}=\{P(i,j), \,P(\Psi(b,j)), \,M(b,j), \,L(i), \,B(i,j), \,A(i;b, j), \,D, \,Q(b,i)\}$$
such that,
\begin{itemize}
\item
$D=(\underline{X_{1}^{2}}-X_{2}X_{0})\Phi(1,2)
-(X_{1}X_{2}-X_{3}X_{0})\Phi(1,1)$\,;
\medskip

\item $B(i,j)=\begin{cases}
(\underline{X_{i}X_{j}}\, -
X_{\epsilon(i,j)}X_{i+j-\epsilon(i,j)})\Psi_{b,3-b}
-(X_{3}^{a}X_{b}-X_{0}^{a+d+1})\Phi(i,j)\\
\hspace*{1.50in} {\rm if} \quad i,j \in [1,2] \,,\\[2mm]
0 \hspace*{1in} ; \hspace*{0.40in} {\rm otherwise}\,;
\end{cases}$
\medskip

\item $A(i;b, j)=\begin{cases}
\underline{X_{i}\Psi(b,j)}-X_{b+i+j-\epsilon(i,b+j)}\Psi(b,\epsilon(i,b+j)-b)
-X_{3}^{a}\Phi(i,b+j)\\
+X_{0}^{a+d}[\Phi(i,j)-\Phi(b+i+j-3,3-b)]\\
\hspace*{1in} {\rm if} \quad i \in [1,3], \,j \in [0,2-b] \quad {\rm and} \quad b\neq 3\,,\\[2mm]
0 \hspace*{1in} ; \hspace*{0.40in} {\rm otherwise}\,;
\end{cases}$
\medskip

\item $ \,L(i)=\begin{cases}
\underline{X_{i}\Phi(2\,,\,2)}\, - \, X_{i+1}\Phi(1\,,\,2)
\,+\,X_{i+2}\Phi(1\,,\,1) \,\, \quad ; \,\, \quad
{\rm if} \quad i, \in [0,1] \,,\\[2mm]
0 \hspace*{2.90in} ; \hspace*{0.20in} {\rm otherwise}\,;
\end{cases}$
\medskip

\item $Q(b,i)=\begin{cases}
\underline{\Psi(1,0)\Phi(2\,,\,2)}+\Psi(1,2)\Phi(1\,,\,1)-
\Psi(1,1)\Phi(1\,,\,2) \quad {\rm if} \quad b=1 \quad {\rm and}
\quad i=1\,,\\
\underline{\Psi(1,1)^{2}}-\Psi(1,2)\Psi(1,0)-X_{3}^{a-1}\Psi(1,2)\Phi(2\,,\,2)
+X_{0}^{a+d-1}\Psi(1,0)\Phi(1\,,\,1)\\[2mm]
\hspace*{0.35in}
-X_{3}^{a-1}X_{0}^{a+d-1}(\Phi(1\,,\,2)^{2}-\Phi(2\,,\,2)\Phi(1\,,\,1))
\quad {\rm if} \quad b=1 \quad {\rm and} \quad i=2\,,\\
\underline{\Psi(2,0)^{2}\Phi(2\,,\,2)}-X_{3}^{a-1}\Psi(2,1)\Phi(2\,,\,2)^{2}
-\Psi(2,1)\Psi(2,0)\Phi(1,2)-X_{3}^{a-1}X_{0}^{a+d-1}\Phi(1\,,\,2)^{3}\\[2mm]
+\Psi(2,1)^{2}\Phi(1\,,\,1)+X_{3}^{a-1}X_{0}^{a+h-1}\Phi(2,2)\Phi(1\,,\,2)\Phi(1\,,\,1)
+X_{0}^{a+d-1}\Psi(2,0)\Phi(1\,,\,1)^{2}\\
\quad \quad \quad {\rm if} \quad b=2
\quad {\rm and} \quad i=1\,,\\
0 \quad \quad \quad \quad \quad \quad \quad {\rm otherwise}\,;\\
\end{cases}$
\medskip

\item $M(b,i)=\begin{cases}
\underline{tX_{0}^{a+d+1}\Psi(b,i)}+\Psi(b,0)\Psi(b,i)
-tX_{3}^{a-1}X_{1+b+i}X_{b-1}\Psi(b,3-b)\\-tX_{3}^{2a}\Phi(b,b+i)
+(-1)^{i+1}tX_{3}^{a-1}X_{0}^{a+d}X_{3b+3i-3}\Phi(3-b-i,3-b-i)\\
\hspace*{0.60in}{\rm if} \quad i \in [0,2-b] \quad {\rm and} \quad b\neq 3\,,\\
0 \hspace*{0.40in} , \hspace*{0.30in}{\rm otherwise}. \\
\end{cases}
$\\
\end{itemize}
}
\end{theorem}

\noindent For our convenience let us set the following:
\begin{enumerate}
\item $V(i,j;q)=\begin{cases}B(i,j) & {\rm if} \quad q=1\,,\\
P(i,j) & {\rm if} \quad q=2\,;\\
\end{cases}$

\item $U(q)=\begin{cases}\Psi(b,3-b) & {\rm if} \quad q=1\,,\\
t & {\rm if} \quad q=2\,;\\
\end{cases}$

\item $u(q)=\begin{cases}\psi(b,3-b) & {\rm if} \quad q=1\,,\\
1 & {\rm if} \quad q=2.\\
\end{cases}$

\item $X_{i}=0 \quad {\rm if} \quad i \notin [0,3]$;

\item $\Phi(i,j)=\Phi(j,i)$;

\item $\Phi(i,j)=0 \quad {\rm if} \quad i,j \notin
[1,2]$;

\item $\Psi(b,j)=0 \quad {\rm
if} \quad j \notin [0,3-b]$;

\item $\phi(i,j)=\phi(j,i)$ \, and \, $V(i,j;q)=V(j,i;q)$.
\end{enumerate}

\noindent The following Lemma will be used for proving Theorem 5.1.

\begin{lemma}
{\it Given $b\in \{1, 2 , 3\}$, let $\mathfrak{Q}_{b}$ be the ideal
in $S$, generated by $\{P(i,j),P(\Psi(b,3-b))\}$\,. A
Gr\"{o}bner Basis for $\mathfrak{Q}_{b}$ is the set
$\widehat{\mathfrak{Q}_{b}}=\{P(i,j), P(\Psi(b,3-b)),L(i),B(i,j),D\}$.}
\end{lemma}

\proof We apply the Buchberger's criterion and show that all the
$S$-polynomials reduce to zero modulo $\widehat{\mathfrak{Q}_{b}}$.
If $\gcd(\,\lm(f), \,\lm(g)\,) = 1$, then the $S$-polynomials reduce
to $0$ modulo $\widehat{\mathfrak{Q}_{b}}$. Let us
consider the other cases, that is when the gcd is not one.
\begin{enumerate}
\item $S(V(1,i;q),L(1))=\Phi(2,2)V(1,i;q)-X_{i}\Psi(b,0)L(1)\,, \quad
{\rm where} \quad i\,\in\,[1,2]\\[2mm]
= -U(q)[\underline{X_{1+i}X_{0}\Phi(2,2)}
-X_{i}X_{2}\Phi(1,2)+X_{i}X_{3}\Phi(1,1)]
-u(q)\Phi(2,2)\Phi(1,i)\\[2mm]
=-X_{1+i}U(q)L(0)+\Phi(1,i)V(2,2;q)$\\[2mm]

\item $S(V(1,i;q),D)=
X_{1}^{i-1}\Phi(1,2)V(1,i;q)- X_{2}^{i-1}U(q)D\,, \quad
{\rm where} \quad i\,\in\,[1,2]\\[2mm]
=-X_{0}U(q)\Phi(1,2)[X_{1}^{i-1}X_{1+i}-X_{2}^{i-1}X_{2}] -
u(q)X_{1}^{i-1}\Phi(1,2)\Phi(1,i)+\phi(1,2)X_{2}^{i-1}U(q)\Phi(1,1)\\[2mm]
=X_{0}^{i-1}\Phi(1,i)V(i,2;q)+X_{2}\Phi(1,i-1)V(1,2;q)+
u(q)\Phi(1,2)L(i-2)$\\[2mm]
${\rm Note\,\, that\, the\, LT\,\, is}\,\,
X_{2}^{i-1}U(q)X_{1}X_{2}\Phi(1,1)\,\, {\rm if}\,\, i=1\,, \, {\rm
and\, the\, LT\,\, is}\,\, X_{2}^{i-1}U(q)X_{2}X_{0}\Phi(1,2)\,\, {\rm
if}\,\, i=2$.\\[2mm]

\item $S(V(1,1;q),V(2,2;q))= X_{2}^{2}V(1,1;q) -
X_{1}^{2}V(2,2;q)\\[2mm]
=-U(q)[X_{2}^{3}X_{0} - \underline{X_{1}^{3}X_{3}}]
+u(q)[X_{1}^{2}\Phi(2,2)-X_{2}^{2}\Phi(1,1)]\\[2mm]
=X_{1} X_{3}V(1,1;q)-X_{2}X_{0}V(2,2;q)+u(q)[X_{1}L(1) - X_{2}L(0)]$\\[2mm]

\item $S(V(1,2;q),V(i,i;q))=X_{i}V(1,2;q)-X_{j}V(i,i;q)\,,
\quad {\rm where} \quad i\in [1,2]\quad {\rm and} \quad j\, \in\, \{1,2\} \setminus \{i\}\\
=-U(q)[X_{i}X_{3}X_{0}-
\underline{X_{j}X_{\epsilon(i,i)}X_{2i-\epsilon(i,i)}}]
-u(q)[X_{i}\Phi(1,2)-X_{j}\Phi(i,i)] \\[2mm]
=X_{3i-3}V(j,j;q) +u(q)L(2-j)$\\[2mm]

\item $S(L(0),L(1))\,=X_{1}L(0)-X_{0}L(1)=-[
\underline{X_{1}^{2}}-X_{2}X_{0}]\Phi(1,2)+\phi(1,2)\Phi(1,1) =-D$\\[2mm]

\item $S(L(1),D)=
X_{1}\Phi(1,2)L(1)-\Phi(2,2)D\\[2mm]
=-\Phi(1,2)[X_{1}X_{2}\Phi(1,2)-X_{1}X_{3}\Phi(1,1)
-\underline{X_{2}X_{0}\Phi(2,2)}]
+\phi(1,2)\Phi(2,2)\Phi(1,1)\\[2mm]
=X_{2}\Phi(1,2)L(0)-\Phi(1,1)[X_{3}L(0)-X_{2}L(1)]$\\[2mm]

\item $S(P(i,j),B(i,j))=\Psi(b,3-b)P(i,j)-tB(i,j)\,,
\quad {\rm where} \quad i,j \in [1,2]\\[2mm]
= -\Psi(b,3-b)\Phi(i,j)+\underline{t\psi(b,3-b)\Phi(i,j)}
=\Phi(i,j)P(\Psi(b,3-b))$\\[2mm]

\item $S(P(i,j),B(l,j))=X_{l}\Psi(b,3-b)P(i,j)-tX_{i}B(l,j)\,,
\quad {\rm where} \quad i\,,\,j\,,\,l \,\in \,[1,2] \, \quad {\rm with} \quad l\neq i\\[2mm]
=-\Psi(b,3-b)[tX_{l}X_{\epsilon(i,j)}X_{i+j-\epsilon(i,j)}+X_{l}\Phi(i,j)
-tX_{i}X_{\epsilon(l,j)}X_{l+j-\epsilon(l,j)}]+tX_{i}\psi(b,3-b)\Phi(l,j)\\[2mm]
=
+X_{i}\Phi(l,j)P(\Psi(b,3-b))+(-1)^{i+j+1}\Psi(b,3-b)[X_{3j-3}P(3-j,3-j)+L(j-1)]$\\[2mm]
${\rm Note\,\, that\, the\, LT\,\, is}\quad
-tX_{l}X_{\epsilon(i,j)}X_{i+j-\epsilon(i,j)}\Psi(b,3-b)\quad  {\rm
if}\quad  l\neq j$\\[2mm]
${\rm and\,the \, LT\,\, is}\quad
tX_{i}X_{\epsilon(l,j)}X_{l+j-\epsilon(l,j)}\Psi(b,3-b)\quad  {\rm
if}\quad  l=j$\\[2mm]

\item $S(P(i,j),P(\Psi(b,3-b)))=X_{3}^{a+1}P(i,j)-X_{i}X_{j}P(\Psi(b,3-b)\,,
\quad {\rm where} \quad i\,,\,j\,\in\,[1,2]\\[2mm]
=-X_{3}^{a+1}[\underline{tX_{\epsilon(i,j)}X_{i+j-\epsilon(i,j)}}+\Phi(i,j)]
+X_{i}X_{j}[tX_{3-b}X_{0}^{a+d}+\Psi(b,3-b)]\\[2mm]
=-X_{\epsilon(i,j)}X_{i+j-\epsilon(i,j)}P(\Psi(b,3-b))+B(i,j)
+X_{3-b}X_{0}^{a+d}P(i,j)$\\[2mm]
\end{enumerate}

\noindent Hence the proof.\qed

\begin{lemma}{\it A minimal Gr\"{o}bner basis for the ideal
$\mathfrak{q}_{b}=\mathfrak{Q}_{b}\cap R_{b}$ is the set $\widehat{\mathfrak{q}_{b}} = \{L(i), B(i,j)\}$.}
\end{lemma}

\proof By the Elimination theorem, a Gr\"{o}bner basis for the ideal $\mathfrak{q}_{b}$ is the
set $\{L(i), B(i,j), D\}$, which contains only those elements of $\widehat{\mathfrak{Q}_{b}}$,
which do not involve the variable $t$. Now, $D=X_{0}L(1)-X_{1}L(0)$ and the leading monomials
of $L(i)$ or $B(i,j)$ do not divide each other. Therefore, by removing $D$ from the above list
we obtain a minimal Gr\"{o}bner basis $\widehat{\mathfrak{q}_{b}} = \{L(i), B(i,j)\}$ for the ideal
$\mathfrak{q}_{b}$.\qed

\begin{corollary}{\it A minimal Gr\"{o}bner basis for the ideal $E_{3}$ is the set
$\widehat{E_3} = \{L(i), B(i,j)\}$.}
\end{corollary}

\proof Note that $\mathfrak{q}_{3}=E_{3}$. Hence, the proof follows
from Lemma 5.3\,.\qed

\begin{remark}
Note that, for $b=3$, the ideal $\wp$ is a prime ideal with
$\mu({\wp})=4=1+{\rm ht}(\wp)$, and therefore an ideal of linear
type by Huneke (1981) and Valla (1980, 1980/81). It is interesting to note that
$\mu({\mathfrak{q}_{b}})=4=1+{\rm ht}(\mathfrak{q}_{b})$, and what we have proved above
shows that $\mathfrak{q}_{b}$ is an ideal of linear type for $b\in
[1,3]$\,, but $\mathfrak{q}_{b}$ is not a prime
ideal if $b\neq 3$. This produces a class of non-prime ideals of linear type which
have the property that $\mu(-)= 1+{\rm ht}(-)$.
\end{remark}

\noindent \textbf{Proof of Theorem 5.1.}

\proof We apply the Buchberger's criterion and show that all the
$S$-polynomials reduce to zero modulo $\widehat{\mathfrak{a}_{b}}$.
If $\gcd(\,\lm(f), \,\lm(g)\,) = 1$, then the $S$-polynomials reduce
to $0$ modulo $\widehat{\mathfrak{a}_{b}}$. Let us
consider the other cases, that is when the gcd is not one.
\medskip

Note that, by Lemma 5.2, every non-zero polynomial $\textbf{H} \in K[t,X_{i},\Psi(b,3-b),\Phi(i,j)]
\subseteq R_{b}[t]$, with $\overline{\varphi_{b}}(\textbf{H})=0$, can be expressed as $\textbf{H}=\sum_{i}
c_{i}H_{i}$, with $c_{i} \in R$, $H_{i} \in \widehat{\mathfrak{Q}_{b}}$ and ${\rm Lm}(\textbf{H}) \geq
{\rm Lm}(c_{i}H_{i})$, whenever $c_{i}\neq 0$. Henceforth, the symbols $\textbf{G}$ and $\textbf{H}$ will only
denote polynomials in $K[t,X_{i},\Psi(b,3-b),\Phi(i,j)] \subseteq R_{b}[t]$, such that
$\overline{\varphi_{b}}(\textbf{H})=0$\,. We use this observation below to prove that the $S$-polynomials converge
to zero. We only indicate the proof for the $S$-polynomial $S(A(i;b,j),A(l;b,j))$, for all othere cases the proof is similar.

\begin{enumerate}
\item $S(A(i;b,j),A(l;b,j)) = X_{l}A(i;b,j)-X_{i}A(l;b,j)\,,\quad
{\rm with }\quad i<l \quad {\rm and}\quad i,l \in [1,3]\\[2mm]
= -X_{l}[X_{b+i+j-\epsilon(i,b+j)}\Psi(b,\epsilon(i,b+j)-b)
+X_{3}^{a}\Phi(i,b+j)
-X_{0}^{a+d}\{\Phi(i,j)-\Phi(b+i+j-3,3-b)\}]\\[2mm]
+X_{i}[X_{b+l+j-\epsilon(l,b+j)}\Psi(b,\epsilon(l,b+j)-b)
+X_{3}^{a}\Phi(l,b+j)
-X_{0}^{a+d}\{\Phi(l,j)-\Phi(b+l+j-3,3-b)\}]$\\[2mm]
= $-\underline{X_{l}X_{b+j+i-\epsilon(i,b+j)}\Psi(b\,,\,\epsilon(i,b+j)-b)}
+X_{i}X_{b+j+l-\epsilon(l,b+j)}\Psi(b\,,\,\epsilon(l,b+j)-b) + \textbf{G}$\\[2mm]
where \,$\textbf{G}$\, is an element of \,$R[t,X_{1},X_{2},X_{3},X_{0},\Psi(b,3-b),\Phi(i,j)]$.
Note that the only monomial of \,$X_{l}A(i;b,j)-X_{i}A(l;b,j)$, which does not belong to
\,$R[t,X_{1},X_{2},X_{3},X_{0},\Psi(b,3-b),\Phi(i,j)]$\, is $-X_{l}X_{b+i+j-\epsilon(i,b+j)}\Psi(b,\epsilon(i,b+j)-b)$ \, if
\, $(i;b,j)=(1;1,0)$. Therefore, every monomial of
\begin{eqnarray*}
\textbf{H} & = & S(A(i;b,j),A(l;b,j))+X_{b+j+i-\epsilon(i,b+j)}A(l;b,\epsilon(i,b+j)-b)\\
           & = & X_{l}A(i;b,j)-X_{i}A(l;b,j)+X_{b+j+i-\epsilon(i,b+j)}A(l;b,\epsilon(i,b+j)-b)
\end{eqnarray*}
belongs to \,$R[t,X_{1},X_{2},X_{3},X_{0},\Psi(b,3-b),\Phi(i,j)]$\,, since every monomial of
\,$A(l;b,\epsilon(i,b+j)-b)$ belongs to \,$R[t,X_{1},X_{2},X_{3},X_{0},\Psi(b,3-b)]$, except
\,$X_{l}\Psi(b,\epsilon(i,b+j)-b)$\, with $(i;b,j)=(1;1,0)$.
\medskip

Moreover, $X_{l}A(i;b,j)$, \,$X_{i}A(l;b,j)$ and \,$X_{b+j+i-\epsilon(i,b+j)}A(l;b,\epsilon(i,b+j)-b)$
belong to ${\rm ker}(\overline{\varphi_{b}})$.
Hence, $\overline{\varphi_{b}}(\textbf{H})=0$. Therefore, we can write
\medskip

$S(A(i;b,j),A(l;b,j)) = X_{l}A(i;b,j)-X_{i}A(l;b,j)\,, {\rm with} \,i<l\, {\rm and} \,i,l \in [1,3];\\[2mm]
= -\underline{X_{l}X_{b+j+i-\epsilon(i,b+j)}\Psi(b\,,\,\epsilon(i,b+j)-b)}
+X_{i}X_{b+j+l-\epsilon(l,b+j)}\Psi(b\,,\,\epsilon(l,b+j)-b) + \textbf{G}\\[2mm]
= -X_{b+j+i-\epsilon(i,b+j)}A(l;b,\epsilon(i,b+j)-b)+\textbf{H}.$
\medskip

\noindent Now one can apply Lemma 5.3 to conclude that there exist \,$c_{i} \in R$ and $H_{i} \in \widehat{\mathfrak{q}_{b}}$, such that
\,$\textbf{H}=\sum_{i}c_{i}H_{i}$.\\[2mm]

\item $S(D,A(1;b,j))=\Psi(b,j)D-X_{1}\Phi(1,2)A(1;b,j)\\[2mm]
=-\underline{X_{2}X_{0}\Psi(b,j)\Phi(1,2)}-\phi(1,2)\Psi(b,j)\Phi(1,1)
+X_{1}X_{0}\Psi(b,1+j)\Phi(1,2)+\textbf{G}\\[2mm]
=-[X_{0}\Phi(1,2)+X_{1}\Phi(1,1)]A(2;b,j)+X_{0}\Phi(1,1)A(3;b,j)
+X_{0}\Phi(1,2)A(1;b,1+j)+\textbf{H}$\\[2mm]

\item $S(L(1),A(1;b,j))=\Psi(b,j)L(1)-\Phi(2,2)A(1;b,j)\\[2mm]
=-\underline{X_{2}\Psi(b,j)\Phi(1,2)}+X_{3}\Psi(b,j)\Phi(1,1)
+X_{0}\Psi(b,1+j)\Phi(2,2)+\textbf{G}\\[2mm]
=-\Phi(1,2)[A(2;b,j)-A(1;b,1+j)]+\Phi(1,1)[A(3;b,j)-A(2;b,1+j)]
+\Psi(b,1+j)L(0) +\textbf{H}$\\[2mm]

\item $S(P(i,l),P(\Psi(b,j))=X_{3}^{a}X_{b+j}P(i,l)-X_{i}X_{l}P(\Psi(b,j))
\quad {\rm where} \quad b+j \notin \{i,l\}\quad {\rm and} \quad i\leq l\\[2mm]
=-\underline{tX_{b+j}X_{\epsilon(i,l)}X_{i+l-\epsilon(i,l)}X_{3}^{a}}-X_{b+j}X_{3}^{a}\Phi(i,l)
+tX_{i}X_{l}X_{j}X_{0}^{a+d}+X_{i}X_{l}\Psi(b,j)\\[2mm]
=-X_{\epsilon(i,l)}X_{i+l-\epsilon(i,l)}P(\Psi(b,j))
-X_{i+l-\epsilon(i,l)}A(\epsilon(i,l);b,j)
+X_{j}X_{0}^{a+d}P(i,l)+X_{i}A(l;b,j)
+\textbf{H}$\\[2mm]

\item  $\displaystyle{S(P(\Psi(b,j)),P(b+j,l))
=X_{l}P(\Psi(b,j))-X_{3}^{a}P(b+j,l)  }$ $\displaystyle{=
\underline{tX_{\epsilon(b+j,l)}X_{b+j+l-\epsilon(b+j,l)}X_{3}^{a}}
-X_{l}\Psi(b,j)+\textbf{G}}$\\[2mm]
$\displaystyle{=-A(l;b,j)-X_{0}P(\Psi(b,3b+5j+3l-5)) +\textbf{H} }$\\[2mm]

\item  $\displaystyle{S(P(\Psi(b,j)),B(b+j,l))
=X_{l}\Psi(b,3-b)P(\Psi(b,j))-tX_{3}^{a}B(b+j,l)  }$\\[2mm]
$\displaystyle{=
\underline{tX_{\epsilon(b+j,l)}X_{b+j+l-\epsilon(b+j,l)}X_{3}^{a}\Psi(b,3-b)}
-X_{l}\Psi(b,j)\Psi(b,3-b)+\textbf{G}}$\\[2mm]
$\displaystyle{=-\Psi(b,3-b)A(l;b,j)-X_{0}\Psi(b,2)P(\Psi(b,3b+5j+3l-5)) +\textbf{H} }$\\[2mm]

\item $S(\,P(\Psi(b,i)\,,\,P(\Psi(b,j))\,)=X_{b+j}P(\Psi(b,i)\,-\,X_{b+i}P(\Psi(b,j))
\quad {\rm assume}\quad i\,<\,j \\[2mm]
=-tX_{i}X_{b+j}X_{0}^{a+d}+\underline{tX_{j}X_{b+i}X_{0}^{a+d}}-X_{b+j}\Psi(b,i)+X_{b+i}\Psi(b,j)\\[2mm]
=-X_{0}^{a+d}[P(i,b+j)-P(j,b +i)]-A(b+j;b,i)+A(b+i;b,j) +\textbf{H}$\\[2mm]

\item $\displaystyle{S(A(i;b,j),P(i,l))=tX_{l}A(i;b,j)-\Psi(b,j)P(i,l)}$\\[2mm]
$\displaystyle{=-tX_{l}X_{b+j+i-\epsilon(i,b+j)}\Psi(b,\epsilon(i,b+j)-b)
+\underline{tX_{\epsilon(i,l)}X_{i+l-\epsilon(i,l)}\Psi(b,j)}
+\Psi(b,j)\Phi(i,l)+\textbf{G}}$\\[2mm]
$\displaystyle{=-tX_{0}A(l;b,3b+5j+3i-5)
+tX_{i+l-\epsilon(i,l)}A(\epsilon(i,l);b,j) +\Phi(i,l)P(\Psi(b,j))
+\textbf{H} }$\\[2mm]

\item $S(A(i;b,j),B(i,l))=X_{l}\Psi(b,3-b)A(i;b,j)-\Psi(b,j)B(i,l)\\[2mm]
=-X_{l}X_{b+j+i-\epsilon(i,b+j)}\Psi(b,\epsilon(i,b+j)-b)\Psi(b,3-b)
+\underline{X_{\epsilon(i,l)}X_{i+l-\epsilon(i,l)}\Psi(b,j)\Psi(b,3-b)}\\[2mm]
+\psi(b,3-b)\Psi(b,j)\Phi(i,l)+\textbf{G}\\[2mm]
=-X_{0}\Psi(b,2)A(l;b,3b+5j+3i-5)
+X_{i+l-\epsilon(i,l)}\Psi(b,3-b)A(\epsilon(i,l);b,j)\\[2mm]
+\Phi(i,l)\Big{[}X_{3}^{a}A(3;b,j)-X_{0}^{a+d}A(3-b;b,j)\Big{]}
+\textbf{H}$\\[2mm]

\item $\displaystyle{S(M(b,j),P(\Psi(b,i)))=X_{3}^{a}X_{b+i}M(b,j)
-X_{0}^{a+d+1}\Psi(b,j)P(\Psi(b,i))}$\\[2mm]
$\displaystyle{=X_{0}^{a+d+1}\Psi(b,j)\Big{[}
\underline{tX_{i}X_{0}^{a+d}}+\Psi(b,i)\Big{]}+X_{3}^{a}X_{b+i}\Psi(b,0)\Psi(b,j) +\textbf{G}}$\\[2mm]
$\displaystyle{=X_{0}^{a+d+1}M(b,i+j)+tX_{0}^{2a+2d+1}A(i;b,j)}$
$\displaystyle{-\Psi(b,3-b)P(\Psi(b,j))\Big{\{}X_{i-2}^{a+d+1}+X_{b-2}^{a+d+1}\Big{\}}}$\\[2mm]
$\displaystyle{+X_{4i-4}^{a+d+1}\Big{[}Q(b,2j)
-\Psi(b,2)P(\Psi(b,2j-2))-X_{0}^{a+d+1}\Phi(2b-1,2b-1)P(\Psi(b,2j-2))\Big{]}}$\\[2mm]
$\displaystyle{+X_{3}^{a-1}X_{b+i}\Psi(b,j)A(3;b,0)
-X_{3}^{a-1}X_{b}\Big{[}X_{b}\Psi(b,3-b)+X_{0}^{a+d}\Phi(b,3-b)\Big{]}P(\Psi(b,j))
+\textbf{H}}$\\[2mm]

\item $S(M(b,j),P(i,l))=X_{i}X_{l}M(b,j)
-X_{0}^{a+d+1}\Psi(b,j)P(i,l)\quad {\rm assume\,\,that}\,\,i \leq l\\[2mm]
=X_{i}X_{l}\Psi(b,0)\Psi(b,j)+X_{0}^{a+d+1}[\underline{tX_{\epsilon(i,l)}X_{i+l-\epsilon(i,l)}}
+\Phi(i,l)]\Psi(b,j)+\textbf{G}\\[2mm]
=X_{i}\Psi(b,0)A(l;b,j)
-X_{i}[X_{b+l+j-\epsilon(l,b+j)}\Psi(b,\epsilon(l,b+j)-b)
+X_{3}^{a}\Phi(l,b+j)]P(\Psi(b,0))\\[2mm]
+X_{i}X_{0}^{a+d}\{\Phi(l,j)-\Phi(b+l+j-3,3-b)\}P(\Psi(b,0))\\[2mm]
+X_{0}^{a+d+1}[
tX_{i+l-\epsilon(i,l)}A(\epsilon(i,l);b,j)-\Phi(i,l)P(\Psi(b,j))]\\[2mm]
+tX_{0}\psi(b,0)A(3j+1;b,7-2b-2j-2l)+\textbf{H}$\\[2mm]

\item $S(L(0),M(b,j))=tX_{0}^{a+d}\Psi(b,j)L(0)-\Phi(2,2)M(b,j)\\[2mm]
=-tX_{0}^{a+d}[\underline{X_{1}\Phi(1,2)}-X_{2}\Phi(1,1)]\Psi(b,j)-\Psi(b,0)\Psi(b,j)\Phi(2,2)
+\textbf{G}\\[2mm]
=-tX_{0}^{a+d}[\Phi(1,2)A(1;b,j)-\Phi(1,1)A(2;b,j)]-
\Psi(b,j)Q(b,b)-Q(b,b-1)\\[2mm]
-\Phi(1,2)Q(b,2j)-M(b,j+1)\Phi(1,2)
+(-1)^{b}\Psi(b,3-b)\Phi(1,b)P(\Psi(b,j))\\[2mm]
+\Psi(1,2)\Phi(1,2)P(\Psi(b,j-1))
-X_{0}^{a+d-1}\Phi(1,j+1)\Phi(1,1)P(\Psi(b,b+j-2))+\textbf{H}$\\[2mm]

\item $S(P(\Psi(b,j)),A(b+j;b,l))=\Psi(b,l)P(\Psi(b,j))-tX_{3}^{a}A(b+j;b,l)
\quad {\rm when} \quad b+j\neq 3\\[2mm]
=-\Psi(b,l)[\underline{tX_{j}X_{0}^{a+d}}+\Psi(b,j)]
+tX_{3}^{a}X_{2b+j+l-\epsilon(b+j,b+l)}\Psi(b,\epsilon(b+j,b+l)-b)+\textbf{G}\\[2mm]
=-tX_{0}^{a+d}A(j;b,l)-M(b,j+l)-Q(b,4l+2j-4)+tX_{3}^{a-1}X_{0}A(j+3;b,b+l)\\[2mm]
+\Psi(b,2)P(\Psi(b,l+j-2))+\Psi(b,1)P(\Psi(b,4b+j-9))
-X_{0}^{a+d-1}\Phi(1,1)P(\Psi(b,4l+3j-7))+\textbf{H}$\\[2mm]

\item $S(P(\Psi(b,j)),A(3;b,l))=\Psi(b,l)P(\Psi(b,j))-tX_{b+j}X_{3}^{a-1}A(3;b,l)
\\[2mm]
=-\Psi(b,l)[\underline{tX_{j}X_{0}^{a+d}}+\Psi(b,j)]
+tX_{b+j}X_{b+l}X_{3}^{a-1}\Psi(b,3-b)+\textbf{G}\\[2mm]
=-tX_{0}^{a+d}A(j;b,l)-M(b,j+l)-Q(b,4l+2j-4)\\[2mm]
+\Psi(b,2)P(\Psi(b,l+j-2))+\Psi(b,1)P(\Psi(b,4b+j-9))
-X_{0}^{a+d-1}\Phi(1,1)P(\Psi(b,4l+3j-7))+\textbf{H}$\\[2mm]

\item $S(M(b,j),A(i;b,j))=X_{i}M(b,j)-tX_{0}^{a+d+1}A(i;b,j)\\[2mm]
=X_{i}[\Psi(b,0)\Psi(b,j)-tX_{b-1}X_{b+j+1}X_{3}^{a-1}\Psi(b,3-b)]
+tX_{b+i+j-\epsilon(i,b+j)}X_{0}^{a+d+1}\Psi(b,\epsilon(i,b+j)-b)+\textbf{G}\\[2mm]
=\Psi(b,0)A(i;b,j)+X_{0}M(b,b+i+j-1)\\[2mm]
-X_{0}\Psi(b,3-b) P(\Psi(b,b+i+j-3) +\Psi(b,3-b)
A(b+i+j-3;b,0)\\[2mm]
-[X_{3}^{a}\Phi(i,b+j)-X_{0}^{a+d}
\{\Phi(i,j)-\Phi(b+i+j-3,3-b)\}]P(\Psi(b,0))+\textbf{H}$\\[2mm]
${\rm
Note\,\,that\, the \,LT\,\,is}\,\,-tX_{i}X_{b-1}X_{b+j+1}X_{3}^{a-1}\Psi(b,3-b)\,\,{\rm
if}\,\,(i;b,j)\neq (1;1,0),\\[2mm]
{\rm and\, the\, LT\,\, is}\,\,
tX_{b+i+j-\epsilon(i,b+j)}X_{0}^{a+d+1}\Psi(b,\epsilon(i,b+j)-b)\,\,{\rm if}\,\,(i;b,j)=(1;1,0)$.\\[2mm]
\end{enumerate}

\noindent Rest of the $S$-polynomial computations and their
reductions modulo $\widehat{\mathfrak{a}_{b}}$ is divided into two cases,
depending on $b=1$ and $b=2$\,.

\noindent \textbf{\underline{Case (i): $b=1$}}

\begin{enumerate}
\item $S(A(i;1,0),A(i;1,1))=\Psi(1,1)A(i;1,0)-\Psi(1,0)A(i;1,1)\\[2mm]
=-X_{1+i-\epsilon(1,i)}\Psi(1,\epsilon(1,i)-1)\Psi(1,1)+X_{i-1}\Psi(1,0)\Psi(1,2)\\[2mm]
-X_{3}^{a}[\Psi(1,1)\Phi(i,1)-\Psi(1,0)\Phi(i,2)]
-X_{0}^{a+d}[\Psi(1,1)\Phi(i-2,2)+\Psi(1,0)\{\Phi(i,1)-\Phi(i-1,2)\}]\\[2mm]
=-X_{3}^{a-1}[\Phi(i,1)A(3;1,1)-\Phi(i,2)A(3;1,0)]
-\Psi(1,2)[A(1+i-\epsilon(1,i);1,1)-A(i-1;1,0)]\\[2mm]
+X_{0}^{a+d}Q(1,i-2) -X_{0}Q(1,i+1)+\textbf{H}$\\[2mm]
${\rm Note \,\, that \, the \,LT \,\, is}\quad
-X_{1+i-\epsilon(i,1)}\Psi(1,1)\Psi(1,\epsilon(i,1)-1)\quad  {\rm
if}\,\,
i=1$,\\[2mm]
${\rm and\, the\, LT\,\, is}\quad  X_{i-1}\Psi(1,0)\Psi(1,2)
\quad  {\rm if} \,\, i\neq 1$.\\[2mm]

\item $S(P(\Psi(1,0)),L(1))=\Phi(2,2)P(\Psi(1,0))-tX_{3}^{a}L(1)\\
=-\underline{tX_{0}^{a+d+1}\Phi(2,2)}-\Phi(2,2)\Psi(1,0)+\textbf{G}
=-Q(1,1)+\Phi(1,2)P(\Psi(1,1))+\textbf{H}$\\[2mm]

\item $S(P(\Psi(1,0)),D)=X_{1}\Phi(1,2)P(\Psi(1,0))-tX_{3}^{a}D\\[2mm]
=-\underline{tX_{2}X_{3}^{a}X_{0}\Phi(1,2)}-X_{1}\Psi(1,0)\Phi(1,2)+\textbf{G}
=-\Phi(1,2)A(1;1,0)+X_{0}\Phi(1,2)P(\Psi(1,1))
+\textbf{H}$\\[2mm]

\item $S(Q(1,1),L(i))=X_{i}Q(1,1)-\Psi(1,0)L(i)\\[2mm]
=-\Phi(1,2)[X_{i}\Psi(1,1)-\underline{X_{i+1}\Psi(1,0)}]
-X_{2+i}\Psi(1,0)\Phi(1,1)+\textbf{G}\\[2mm]
=-\Phi(1,2)[A(i;1,1)-A(i+1;1,0)]-\Phi(1,1)A(i+2;1,0)+\textbf{H}$\\[2mm]

\item $S(Q(1,1),A(i;1,0))=X_{i}Q(1,1)-\Phi(2,2)A(i;1,0)\\[2mm]
=-X_{i}\Psi(1,1)\Phi(1,2)+X_{1+i-\epsilon(i,1)}\Psi(1,\epsilon(i,1)-1)\Phi(2,2)+\textbf{G}\\[2mm]
=-\Phi(1,2)[A(i;1,1)-A(1;1,i)]-\Phi(1,1)A(2;1,i)+\Psi(1,i)L(0)+\textbf{H}$\\[2mm]
${\rm Note \,\, that \, the\, LT \,\, is}\quad -X_{i}\Psi(1,1)\Phi(1,2)\quad {\rm if} \quad i \neq 1$,\\[2mm]
${\rm and \, the\, LT \,\, is}\quad
X_{1+i-\epsilon(i,1)}\Psi(1,\epsilon(i,1)-1)\Phi(2,2)
\quad {\rm if} \quad i = 1$.\\[2mm]

\item $S(Q(1,1),M(1,0))=tX_{0}^{a+d+1}Q(1,1)-\Phi(2,2)M(1,0)\\[2mm]
=-\underline{tX_{0}^{a+d+1}\Psi(1,1)\Phi(1,2)}-\Psi^{2}(1,0)\Phi(2,2)+\textbf{G}\\[2mm]
=-\Phi(1,2)M(1,1)-\Psi(1,0)Q(1,1)-\Psi(1,2)\Phi(1,1)P(\Psi(1,0))+\textbf{H}$\\[2mm]

\item $S(Q(1,2),A(i;1,1))=X_{i}Q(1,2)-\Psi(1,1)A(i;1,1)\\[2mm]
=-\Psi(1,2)[\underline{X_{i}\Psi(1,0)}-X_{i-1}\Psi(1,1)]
+X_{3}^{a}\Psi(1,1)\Phi(i,2)
\\[2mm]
+X_{0}^{a+d-1}[X_{i}\Psi(1,0)\Phi(1,1)-X_{0}\{\Phi(i,1)-\Phi(i-1,2)\}\Psi(1,1)]+\textbf{G}\\[2mm]
=-\Psi(1,2)[A(i;1,0)-A(i-1;1,1)]+X_{3}^{a-1}\Phi(i,2)A(3;1,1)\\[2mm]
+X_{0}^{a+d-1}[\Psi(1,1)L(i-3)+\Phi(1,2)A(i-2;1,1)-\Phi(1,1)A(8-2i;1,1)+\textbf{H}$\\[2mm]

\item $S(Q(1,2),M(1,1))=tX_{0}^{a+d+1}Q(1,2)-\Psi(1,1)M(1,1)\\[2mm]
=-tX_{0}^{a+d+1}[\underline{\Psi(1,0)\Psi(1,2)}-X_{0}^{a+d-1}\Psi(1,0)\Phi(1,1)]-\Psi(1,0)\Psi^{2}(1,1)\\[2mm]
+tX_{3}^{a}\Psi(1,1)[X_{0}\Psi(1,2)-X_{0}^{a+d}\Phi(1,1)+X_{3}^{a}\Phi(1,2)]+\textbf{G}\\[2mm]
=-[\Psi(1,2)-X_{0}^{a+d-1}\Phi(1\,,\,1)]M(1,0)
+X_{0}^{a+d-1}X_{3}^{a-1}[ \Phi^{2}(1,2)-\Phi(2,2)\Phi(1,1)]P(\Psi(1,0))\\[2mm]
+X_{3}^{a-1}\Psi(1,2)\Phi(2,2)P(\Psi(1,0)) -\Psi(1,0)Q(1,2)
+tX_{3}^{a-1}[X_{0}\Psi(1,2)+X_{0}^{a+d}\Phi(1,1)+X_{3}^{a}\Phi(1,2)]A(3;1,1)\\[2mm]
-tX_{3}^{2a-1}[\Psi(1,2)L(1)-X_{0}^{a+d-1}\{
\Phi(1,2)L(0)+\Phi(1,1)L(1) \} ]+\textbf{H}$\\[2mm]

\item $S(M(1,0),M(1,1))=\Psi(1,1)M(1,0)-\Psi(1,0)M(1,0)\\[2mm]
=-tX_{3}^{a-1}\Psi(1,1)[ X_{2}X_{0}\Psi(1,2)+X_{3}^{a+1}\Phi(1,1)
+X_{0}^{a+d}X_{0}\Phi(2,2)]\\[2mm]
+tX_{3}^{a-1}\Psi(1,0)[
\underline{X_{3}X_{0}\Psi(1,2)}+X_{3}^{a+1}\Phi(1,2)
-X_{0}^{a+d}X_{3}\Phi(1,1)]+\textbf{G}\\[2mm]
=-tX_{3}^{a-1}X_{0}\Psi(1,2)[A(2;1,1)-A(3;1,0)]
-tX_{3}^{2a-1}[\Phi(1,1)A(3;1,1)-\Phi(1,2)A(3;1,0)]\\[2mm]
-tX_{3}^{a-1}X_{0}^{a+d}[\Psi(1,1)L(0)
+\Phi(1,2)A(1;1,1)-\Phi(1,1)A(2;1,1)+\Phi(1,1)A(3;1,0)]+\textbf{H}$\\[2mm]
\end{enumerate}

\noindent \textbf{\underline{Case(ii): $b=2$}}

\begin{enumerate}
\item $S(M(2,0),Q(2,1))=\Psi(2,0)\Phi(2,2)M(2,0)
-tX_{0}^{a+d+1}Q(2,1)\\[2mm]
=\Psi(2,0)\Phi(2,2)[\Psi(2,0)\Psi(2,0)
-\underline{tX_{3}^{a}X_{1}\Psi(2,1)}-tX_{3}^{2a}\Phi(2,2)
-tX_{3}^{a-1}X_{0}^{a+d}X_{3}\Phi(1,1)]\\[2mm]
+tX_{0}^{a+d+1}[ X_{3}^{a-1}\Psi(2,1)\Phi^{2}(2\,,\,2)
+\Psi(2,0)\Psi(2,1)\Phi(1,2)]
+\textbf{G}\\[2mm]
=\Psi(2,0)Q(2,1)+\Psi(2,1)\Phi(1,2)M(2,0)
-X_{0}^{a+d-1}\Phi^{2}(1,1)M(2,0)-X_{3}^{a-1}\Psi(2,0)\Phi^{2}(2,2)P(\Psi(2,1))\\[2mm]
+\Psi^{2}(2,1)\Phi(1,1)P(\Psi(2,0))
-X_{3}^{a-1}X_{0}^{a+d-1}[\Phi^{3}(1,2)
-\Phi(2,2)\Phi(1\,,\,2)\Phi(1\,,\,1)]P(\Psi(2,0))\\[2mm]
-tX_{3}^{a-1}\Phi(2,2)[X_{1}\Psi(2,1)A(3;2,0)
+X_{0}^{a+d}\{\Phi(1\,,\,1)A(3;2,0) +\Phi(2,2)A(1;2,0)\}]
+\textbf{H}$\\[2mm]

\item $S(Q(2,1),A(i;2,0))=X_{i}Q(2,1)-\Psi(2,0)\Phi(2,2)A(i;2,0)\\[2mm]
=-X_{i}\Psi(2,0)[\Psi(2,1)\Phi(1,2)-X_{0}^{a+d-1}\Phi(1\,,\,1)^{2}]
\\[2mm]
+\Psi(2,0)\Phi(2,2)[\underline{X_{i-1}\Psi(2,1)} +X_{3}^{a}\Phi(i,2)
+X_{0}^{a+d}\Phi(i-1,1)]+\textbf{G}\\[2mm]
=-[\Psi(2,1)\Phi(1,2)-X_{0}^{a+d-1}\Phi^{2}(1,1)]A(i;2,0)
+\Psi(2,1)\Phi(2,2)A(i-1;2,0)\\[2mm]
+X_{3}^{a-1}\Phi(i,2)\Phi(2,2)A(3;2,0)
+X_{0}^{a+d-1}\Phi(i-1,1)\Psi(2,0)L(0)\\[2mm]
+X_{0}^{a+d-1}\Phi(i-1,1)[\Phi(1,2)A(1;2,0)
-\Phi(1,1)A(2;2,0)]+\textbf{H}$\\[2mm]

\item  $S(Q(2,1),L(i))=X_{i}Q(2,1)-\Psi(2,0)^{2}L(i)\\[2mm]
=-X_{i} \Psi(2,0)[\Psi(2,1)\Phi(1,2)-
X_{0}^{a+d-1}\Phi(1\,,\,1)^{2}]
+\Psi^{2}(2,0)[\underline{X_{1+i}\Phi(1,2)}-X_{2+i}\Phi(1,1)]+\textbf{G}\\[2mm]
=\Psi(2,0)[\Phi(1,2)A(i+1;2,0)-\Phi(1,1)A(i+2;2,0)]\\[2mm]
-\Psi(2,1)[\Phi(1,1)A(i+1;2,0)-\Phi(1,2)A(i;2,0)]\\[2mm]
+X_{3}^{a-1}[\Phi(2,1+i)\Phi(1,2)-\Phi(2,2+i)\Phi(1,1)]A(3;2,0)\\[2mm]
+X_{0}^{a+d-1}\Phi(1,1)^{2}A(i;2,0)+\textbf{H}$\\[2mm]
\end{enumerate}

Hence the proof.\qed

\begin{theorem}
{\it Given $b\in \{1, 2\}$, a Gr\"{o}bner Basis for the ideal $E_{b}$
is the set
$$\widehat{E_{b}}=\{A(i;b,j),B(i,j),D,L(i),Q(b,i)\}.$$}
\end{theorem}

\proof Note that, $\widehat{E_{b}}=\widehat{\mathfrak{a}_{b}}\cap R_{b}$\,.\qed
\bigskip

Furthermore, if $b\in \{1, 2\}$, \,$b+l=2$ and \,$i,j \in [1,2]$, we have
$$B(i,j)=X_{i+1}A(j;b,l)-X_{j}A(i+1;b,l)-X_{3}^{a}L(i+j-1)
-X_{0}^{a+d}L(2i+2j-5b)+X_{0}^{a+d}L(7-b-i-j).$$
Therefore, a smaller set $\widehat{\widehat{E_{b}}} =
\{A(i;b,j),L(i),Q(b,i),D\}$ generates the ideal $E_{b}$.

\section{Smoothness of Blowups}

\noindent Let $E$ and $\mathfrak{P}=\langle Y_{1},\ldots,Y_{n-1} \rangle$ be
prime ideals of a ring $N=K[Y_{1},\ldots,Y_{n}]$  and $E \subseteq
\mathfrak{P}$\,. Let $\mathcal{J}_{\mathfrak{P}}$ denote the
Jacobian matrix of the ideal $E$ , taken modulo $\mathfrak{P}$.
Given an indeterminate $\zeta \in \{Y_{1},\ldots,Y_{n}\}$, let
$C_{\zeta}$ denote the column in the matrix
$\mathcal{J}_{\mathfrak{P}}$, corresponding to the indeterminate
$\zeta$. Then, it is obvious from the construction of
$\mathfrak{P}$, that, the column $C_{\zeta}$ is non-zero if and only
if there exists a polynomial $F \in E$ such that $F$ has at least
one term of the form $k\zeta Y_{n}^{l}$, for some $k\in K$ and $l
\in \mathbb{N}$\,.
\medskip

Before we prove our last theorem let us record the following
observations:

\begin{enumerate}
\item
$F\, \in \,\widehat{E_1}$ implies that no term of $F$ is an element
of the set
$$\{X_{2}\Phi(2,2)\,,\,X_{3}\Phi(2,2)\,,\,\Psi(1,1)\Phi(2,2)\,,\,
\Psi(1,2)\Phi(2,2)\,,\,\Phi(1,2)\Phi(2,2)\,,\,\Phi(1,1)\Phi(2,2)\,,\,\Phi(2,2)^{l}\}.$$

\item
$F\,\in\,\widehat{E_2}$  implies that no term of $F$ is an element
of the set
$$\{X_{2}\Phi(2,2)\,,\,X_{3}\Phi(2,2)\,,\,\Psi(2,1)\Phi(2,2)\,,\,
\Psi(2,0)\Phi(2,2)\,,\,\Phi(1,2)\Phi(2,2)\,,\,\Phi(1,1)\Phi(2,2)\,,\,\Phi(2,2)^{l}\}.$$

\item
$F\,\in\,\widehat{E_3}$ implies that no term of $F$ is an element
of the set $$\{X_{1}\Psi(3,0)
,X_{2}\Psi(3,0),X_{3}\Psi(3,0),X_{0}\Psi(3,0)\,,\,
\Psi(2,2)\Psi(3,0)\,,\,\Phi(1,2)\Psi(3,0)\,,\,\Phi(1,1)\Psi(3,0)\,,\,\Psi(3,0)^{l}\}.$$

\end{enumerate}

\begin{theorem}{\it
${\rm Proj}\,\mathcal{R}(\wp)$ is not smooth.}
\end{theorem}

\proof Let us write
\medskip

$\mathfrak{P}_{b} =\begin{cases} \langle X_{1}, X_{2}, X_{3} ,
X_{0}, \Psi(1,0), \Psi(1,1),\Psi(1,2), \Phi(1,2), \Phi(1,1)
\rangle\,, \quad {\rm if} \quad b=1;\\
\langle X_{1}, X_{2}, X_{3} , X_{0}, \Psi(2,0), \Psi(2,1),
\Phi(1,2), \Phi(1,1) \rangle\,, \quad {\rm if} \quad b=2;\\
\langle X_{1}, X_{2}, X_{3} , X_{0}, \Phi(2,2), \Phi(1,2), \Phi(1,1)
\rangle\,, \quad {\rm if} \quad b=3.\\
\end{cases}$
\medskip

\noindent It is clear that $\mathfrak{P}_{b}$ is a homogeneous prime
ideal of $R_{b}$, containing $E_{b}$. Let
$\mathcal{J}_{\mathfrak{P}_{b}}$ denote the Jacobian matrix, taken
modulo $\mathfrak{P}_{b}$. Now we use the preceding observations to conclude that
\medskip

$\bullet$\noindent $C_{\zeta}$ \,is non-zero if and only if \,$\zeta
\, \in \,
\{X_{0},X_{1},\Psi(1,0)\}$, for $b=1$.
\smallskip

$\bullet$\noindent $C_{\zeta}$ \,is non-zero if and only if \,$\zeta
\, \in \,
\{X_{0},X_{1}\}$, for $b=2$.
\smallskip

$\bullet$\noindent $C_{\zeta}$ \,is zero if \,$\zeta \, \in \,
\{X_{1}, X_{2}, X_{3} , X_{0}, \Phi(2,2) , \Phi(1,2), \Phi(1,1)\}$, for $b=3$.
\smallskip

\noindent Hence, the rank of the matrix
$\mathcal{J}_{\mathfrak{P}_{b}}$ is $\begin{cases} 3 \quad {\rm when} \quad b=1;\\
2 \quad {\rm when} \quad b=2;\\
0 \quad {\rm when} \quad b=3;\\
\end{cases}$\\
\noindent and the height of the ideal
$(E_{b})_{(\mathfrak{P}_{b})}$ in the localized ring
$(R_{b})_{(\mathfrak{P}_{b})}$ is $6-b$. Therefore,
$\left(\mathcal{R}(\wp)\right)_{(\mathfrak{P}_{b})}=\left(R_{b}/E_{b}\right)_{(\mathfrak{P}_{b})}$
is not regular by Theorem 3.4. Hence, \,${\rm Proj}
\,\mathcal{R}(\wp)$ \,is not smooth\,.\qed

\end{document}